\documentstyle[12pt]{article}
\begin{document}
\def\beq{\begin{equation}}
\def\eeq{\end{equation}}
\def\bea{\begin{eqnarray}}
\def\eea{\end{eqnarray}}
\def\nn{\nonumber}
\begin{flushright}
q-alg/9504021
\end{flushright}
\begin{center}
{\large \bf FINITE-DIMENSIONAL CALOGERO REPRESENTATION OF THE
$q$-DIFFERENTIAL OPERATOR}
\end{center}

\bigskip

\begin{center}
{R. CHAKRABARTI} \\
{\em Department of Theoretical Physics, University of Madras, \\
Guindy Campus, Madras - 600025, India}

\vspace{.25cm}

{R. JAGANNATHAN}\footnote{E-mail: jagan@imsc.ernet.in} \\
{\em The Institute of Mathematical Sciences, \\
CIT Campus, Tharamani, Madras - 600113, India}
\end{center}

\vspace{2cm}

\baselineskip14pt

\noindent
A finite-dimensional matrix representation of the Jackson $q$-differential
operator $D_q$, defined by $D_qf(x)$ $=$ $(f(qx)-f(x))/(x(q-1))$, is written
down following Calogero.  Such a representation of $D_q$ should have
applications in $q$-analysis leading to the corresponding extensions of the
numerous results of Calogero's work.

\vspace{1cm}

\noindent{\bf 1. Introduction}

\renewcommand{\theequation}{1.{\arabic{equation}}}
\setcounter{equation}{0}

\bigskip

\noindent
Recently Calogero$^1$ has extended to multidimensions his work on a
convenient finite-dimensi-onal matrix representation of the
differential
operator$^{2,3}$.  It may be recalled that the earlier findings in the
one dimensional case
led to several remarkable results related to matrix theory, integrable
dynamical systems, classical polynomials, special functions, numerical
treatment of Sturm-Liouville eigenvalue problems etc. (see Ref. 1 for
detailed bibliography).

The purpose of the present short note is to write down the Calogero
matrix representation of the Jackson $q$-differential operator in the
one dimensional case, defined by
%
%1.1{Dq}
\begin{equation}
D_qf(x) = \frac{f(qx)-f(x)}{x(q-1)} = \left( \frac{q^{x\frac{d}{dx}}
-1}{x(q-1)} \right) f(x)\,.
\label{Dq}
\end{equation}
We do this in view of the central role played by $D_q$ in
$q$-analysis (see, e.g., Refs. 4-6) and the recent interest in
$q$-analysis in relation to the theory of quantum algebras,
possible deformations of the current framework of quantum mechanics,
etc., (see, e.g., Refs. 7-11).  We hope that such a representation of
$D_q$ would lead to extensions of Calogero's remarkable results.
We shall assume $q$ to be generic throughout this note.

\vspace{1cm}

\noindent{\bf 2. Calogero matrices for $x$ and $d/dx$}

\renewcommand{\theequation}{2.{\arabic{equation}}}
\setcounter{equation}{0}

\bigskip

\noindent
Calogero$^1$ gives a straightforward prescription for obtaining the
required convenient finite-dimensional representation of any linear
differential operator ${\cal A}(x,d/dx)$ within a chosen scheme of
interpolation.  Here, we shall adhere to Calogero's original scheme
of the Lagrange interpolation$^{2,3}$.  In this section, we recall
some of the basic formulae of Calogero which will help us write down,
in the next section, the corresponding Calogero matrix representation
for $D_q$.

With $\left\{ f_j \mid j=1,2,\ldots ,n \right\}$ denoting $n$ given
numbers, the `interpolational function'
%
%2.1{fx}
\begin{equation}
f(x) = \sum_{j=1}^{n}\,f_j \delta^{(j)}(x)\,,
\label{fx}
\end{equation}
where $\left\{ \delta^{(j)}(x) \mid j=1,2,\ldots ,n \right\}$
are the `interpolational polynomials'
%
%2.2{ips}
\begin{equation}
\delta^{(j)}(x) = \prod_{k(\neq j)=1}^{n}
\,\frac{x-x_k}{x_j-x_k}\,,
\label{ips}
\end{equation}
is such that $f(x)$ takes the values $\left\{ f_j \mid j=1,2,
\ldots ,n \right\}$, respectively, at the chosen $n$ distinct points
(`nodes') $\left\{ x_j \mid j=1,2,\ldots ,n\right\}$:
%
%2.3{fj}
\begin{equation}
f\left( x_j \right) = f_j\,, \qquad j=1,2,\ldots ,n.
\label{fj}
\end{equation}

The Lagrange interpolational polynomials $\left\{ \delta^{(j)}(x) \mid
j=1,2,\ldots ,n \right\}$ are linear combinations of the monomials
$\left\{ x^{j-1} \mid j=1,2,\ldots ,n \right\}$ (`seeds' of the Lagrange
interpolation), with the structure
%
%2.4{lips}
\begin{equation}
\delta^{(j)}(x) = \sum_{k=1}^{n}\,C_{jk} x^{k-1}\,, \qquad
j=1,2,\ldots ,n\,,
\label{lips}
\end{equation}
where the matrix $C = \left[ C_{jk} \right]$, with $j$ as the
row index and $k$ as the column index, is the inverse of the Vandermonde
matrix $V = \left[ V_{jk} \right] = \left[ x_j^{k-1} \right]$.  Further,
defining
%
%2.5{f,g}
\begin{equation}
\langle f,g \rangle = \sum_{j=1}^{n}\,f\left( x_j \right)g\left( x_j
\right)\,,
\label{f,g}
\end{equation}
one has the orthogonality relation
%
%2.6{orthorel}
\begin{equation}
\left\langle \delta^{(j)} , \delta^{(k)} \right\rangle =
\delta_{jk}\,,
\qquad j,k=1,2,\ldots ,n.
\label{orthorel}
\end{equation}

Now, one has
%
%2.7{deltaj,k}
\begin{equation}
\left\langle \delta^{(j)} ,x\delta^{(k)} \right\rangle = x_j
\delta_{jk}\,, \qquad j,k = 1,2,\ldots ,n.,
\label{deltaj,k}
\end{equation}
or, in other words, the matrix
%
%2.8{X}
\begin{equation}
X = \left[ X_{jk} \right] = \left[ x_j \delta_{jk} \right]
= {\rm diag}\left( x_j \right)\,, \qquad
j,k = 1,2,\ldots ,n.,
\label{X}
\end{equation}
represents the operator `multiplication by $x$' in the
basis $\left\{ \delta^{(j)} \right.$ $\mid $ $j=1,2,$ $\ldots ,$
$\left. n \right\}$.  This
means that, as far as we are concerned with the values of functions
at the distinct chosen nodes $\left\{ x_j \mid j=1,2,\ldots ,n
\right\}$, the interpolational polynomials $\left\{ \delta^{(j)}
\mid j=1,2,\ldots ,n \right\}$ are like the Dirac delta functions
and are `eigenfunctions' of $x$ in a sense.  The corresponding
Calogero matrix for the differential operator $d/dx$ is obtained
from the definition
%
%2.9{D}
\begin{equation}
D_{jk} = \left\langle \delta^{(j)} , \frac{d}{dx} \delta^{(k)}
\right\rangle = \left. \frac{d}{dx} \delta^{(k)} \right| _{x=x_j}\,.
\label{D}
\end{equation}
It should be noted that the set $\left\{ \delta^{(j)}
\mid j=1,2,\ldots ,n \right\}$ is closed under differentiation,
namely, the derivative of any $\delta^{(j)}$ can be expressed as a
linear combination of $\left\{ \delta^{(k)} \mid k=1,2,\ldots ,n
\right\}$ with constant coefficients.  This is a consequence of
(\ref{lips}) and similar closure property of the set of seeds $\left\{
x^{j-1} \mid j=1,2,\ldots n \right\}$ under differentiation.
Explicitly,
%
%2.10{DBZB}
\begin{equation}
D = BZB^{-1}\,,
\label{DBZB}
\end{equation}
where
%
%2.11{B}
\begin{equation}
B = {\rm diag}\left( b_j \right)\,, \qquad
b_j = \prod_{k(\neq j)=1}^{n}\,\left( x_j - x_k \right)\,,
\label{B}
\end{equation}
and
%
%2.12{Z}
\begin{equation}
Z_{jk} = \left\{
\begin{array}{ll}
\left( x_j-x_k \right)^{-1} & {\rm if}\ \ j\neq k\,, \\
                            &        \\
\sum_{k(\neq j)=1}^{n}\,\left( x_j-x_k \right)^{-1} & {\rm if}\ \
j=k\,.
\end{array}
\right.
\label{Z}
\end{equation}
It is clear that for any other differential operator,
${\cal A}\left( x,x/dx \right)$, the corresponding Calogero matrix
will be given by the correspondence rule
%
%2.13{calA}
\begin{equation}
{\cal A}\left( x, \frac{d}{dx} \right) \rightarrow A = {\cal A}
(X,D)\,.
\label{calA}
\end{equation}

It follows from the earlier work of Calogero$^{2,3}$ that another
useful form of the matrix $D$ is:
%
%2.14{D2}
\begin{equation}
D = X^{-1}VNV^{-1}\,,\
\label{D2}
\end{equation}
where $V$ is the Vandermonde matrix and
%
%2.15{N}
\begin{equation}
N = {\rm diag}(j-1)\,, \qquad j=1,2,\ldots ,n\,,
\label{N}
\end{equation}
provided $x_j \neq 0$ for any $j =1,2,\ldots ,n$.

\vspace{1cm}

\noindent{\bf 3. Calogero matrix for $D_q$}

\renewcommand{\theequation}{3.{\arabic{equation}}}
\setcounter{equation}{0}

\bigskip

\noindent
The first observation one can make with reference to the
$q$-differential operator (\ref{Dq}) is that the seeds $\left\{ x^{j-1}
\mid j=1,2,\ldots ,n \right\}$ and hence the polynomials
$\left\{ \delta^{(j)} \mid j=1,2,\ldots ,n \right\}$ are closed
under $q$-differentiation also:
%
%3.1{Dqx}
\begin{equation}
D_qx^{j-1} = [j-1]_q x^{j-2}\,, \qquad j=1,2,\ldots ,
\label{Dqx}
\end{equation}
with the Heine `basic number' $[m]_q$ defined by
%
%3.2{qno}
\begin{equation}
[m]_q = \frac{q^m - 1}{q-1}\,.
\label{qno}
\end{equation}
Hence, it is clear that the Calogero matrix for $D_q$,
say ${\cal D}$, is given by
%
%3.3{calD}
\begin{equation}
{\cal D}_{jk} = \left\langle \delta^{(j)} ,D_q\delta^{(k)}
\right\rangle = \left. D_q\delta^{(k)}(x) \right|_{x=x_j}\,.
\label{calD}\end{equation}
We shall assume that $x_j \neq 0$, for any $j=1,2,\ldots ,
n$.  Then, using the prescription (\ref{calA}), along with
(\ref{X}) and (\ref{D2}),
in (\ref{Dq}), one can see easily that
the explicit form of ${\cal D}$ is as follows:
%
%3.4{Dqrep}
\begin{equation}
{\cal D} = X^{-1}V[N]_qV^{-1}\,, \qquad
[N]_q = {\rm diag}\left( [j-1]_q \right)\,, \ \ \ \ j=1,2,\ldots ,n\,,
\label{Dqrep}
\end{equation}
where $V$ is the Vandermonde matrix.  In the limit $q$
$\rightarrow$ $1$, ${\cal D}$ becomes Calogero's $D$ in (\ref{D2}).  As
is obvious, the matrix ${\cal N}$ $=$ $X{\cal D}$ has the first $n$
basic numbers, corresponding respectively to the first $n$
nonnegative integers, as its eigenvalues independent of the choice
of $\left\{ x_j \mid j=1,2,\ldots ,n \right\}$: this is the
$q$-analogue of Calogero's result that the matrix $XZ$ (or $XD$,
since $Z$ $=$ $B^{-1}DB$ and $B$ commutes with $X$) has the first
$n$ nonnegative integers as its eigenvalues independent of the
choice of $\left\{ x_j \mid j=1,2,\ldots ,n \right\}$.

It may be noted that in $q$-analysis also, the matrix representing
the operator `multiplication by $x$' will be $X$ $=$
${\rm diag}\left( x_j \right)$, $j=1,2,\ldots ,n$.  Only wherever
the differential operator $d/dx$ has to be replaced by $D_q$, the
Calogero matrix $D$ in (\ref{D2}) will have to be replaced by
its $q$-analogue ${\cal D}$ in (\ref{Dqrep}).

\vspace{1cm}

\noindent{\bf 4. Conclusion}

\bigskip

\noindent
To summarize, in this short note we have obtained the finite
dimensional Calogero matrix representation of the $q$-differential
operator $\left( D_q \right)$ that would replace the finite
dimensional Calogero matrix representation of the differential
operator $\left( d/dx \right)$ in studies on the $q$-analogues
of the numerous results of Calogero relating to matrix theory,
integrable dynamical systems, classical polynomials, special
functions, and so on.  The operator `multiplication by $x$' would
have the same Calogero matrix representation, namely, $x$
$\rightarrow$ $X$ $=$ ${\rm diag}\left( x_j \right)$, with
$j=1,2,\ldots ,n$, in the $q$-analysis also (the only restriction
we are required to have is that none of the nodes $\left\{ x_j
\mid j=1,2,\ldots ,n \right\}$ should be at the origin).

\vspace{2cm}

\noindent{\bf Acknowledgments}

\bigskip

\noindent
One us (R.J) wishes to acknowledge gratefully that he got interested
in the Calogero approach to differential calculus mainly from the
discussions, about a decade ago, with Prof. T.S. Santhanam whose
work (see, e.g., Ref. 12) on a formulation (or an approximation) of
quantum mechanics using finite dimensional Hilbert spaces has an
algebraic structure related to the Calogero formalism.  He wishes to
thank Prof. Guido Vanden Berghe for the kind hospitality he enjoyed at
the Department Applied Mathematics and Computer Science, University of
Ghent, where part of this work was done.  It is a pleasure to thank
Prof. Joris Van der Jeugt for fruitful discussions.  He also
acknowledges the support from E.E.C. (contract No. CI1$^*$-CT92-0101)
which enabled him to visit Gent.

\vspace{1.25cm}

\noindent{\bf References}

\begin{enumerate}

%1
\item
F. Calogero, {\em J. Math. Phys.} {\bf 34}, 4704 (1993).

%2
\item
F. Calogero, {\em J. Math. Phys.} {\bf 22}, 919 (1981).

%3
\item
F. Calogero, {\em Rend. Semin. Mat., Univ. Polit. Torino}, October 1984
(Special issue devoted to the Proceedings of the International Conference
on Special Functions: Theory and Computation) p.23.

%4
\item
H. Exton, {\it q-Hypergeometric Functions and Applications} (Ellis Horwood,
Chichester, England, 1983).

%5
\item
G.E. Andrews, {\it q-Series: Their Development and Applications in
Analysis, Number Theory, Combinatorics, Physics and Computer Algebra}
(American Math. Soc., Providence, RI, U.S.A., 1986).

%6
\item
G. Gasper and M. Rahman, {\it Basic Hypergeometric Series}, Encyclopaedia
of Mathematics and its Applications Vol.35 (Cambridge University Press.,
Cambridge, England, 1990).

%7
\item
M. Arik, {\em Z. Phys. C: Particles and Fields} {\bf 51}, 627 (1991).

%8
\item
R. Floreanini and L. Vinet, {\em Ann. Phys.} {\bf 221}, 53 (1993).

%9
\item
R. Floreanini and L. Vinet, {\em J. Phys. A: Math. Gen.} {\bf 26}, L611
(1993).

%10
\item
J.A. Minahan, {\em Mod. Phys. Lett.} {\bf A5}, 2635 (1990).

%11
\item
R. Chakrabarti, R. Jagannathan and R. Vasudevan, {\em Mod. Phys. Lett.}
{\bf A8}, 2695 (1993).

%12
\item
T.S. Santhanam, {\em Found. Phys.} {\bf 7}, 121 (1977).

\end{enumerate}

\noindent
{\em Note}:  After completion of this work, we have learnt that a
closely related discussion of the `finite dimensional representation of
the $Q$ operator' is to be found in a paper by Calogero and Ji Xiaoda
with the title "Solvable (nonrelativistic, classical) $n$-body problems
in multidimensions - II" due to appear in the Proceedings of a Meeting
on Nonlinear Dynamics (Pavullo nel Frignano, Italy, May 1994) Eds.
M. Costato, A. Degasperis and M. Milano (Publ.: Editrice Compositon, the
Publishers of Nuovo Cimento).

\end{document}